%% file: complast.tex
\documentclass[12pt]{amsart}

\usepackage{amsmath,amsfonts,amsthm,amsopn}
\usepackage{graphicx}
\usepackage{epsfig,verbatim}

\newtheorem{Def}{Definition}

\newtheorem{Th}{Theorem}
\newtheorem{Pro}{Proposition}
\input{monimacro}
\begin{document}
\title{Compactness Criteria in Function Spaces}
\author{Monika D\"orfler, Hans G. Feichtinger and Karlheinz Gr\"ochenig}
\address{ Institut f. Mathematik, Universit\"at Wien, Strudlhofg. 4,
  A-1090 Wien,   Austria }
\address{Dept. of Mathematics,
University of Connecticut, Storrs, CT  06269-3009, USA}
\email{monika.doerfler@univie.ac.at, hans.feichtinger@univie.ac.at}
\email{groch@math.uconn.edu}
 \subjclass{46B50,42B35}
\date{}
\keywords{Compactness, Besov spaces,  modulation
  spaces, Bargmann-Fock spaces, coorbit spaces,  wavelet transform,
  short-time Fourier   transform}
\maketitle

\begin{abstract}
 The classical
criterion for compactness in Banach  spaces of functions can be
reformulated into  a simple tightness condition in the
time-frequency domain. This description preserves  more explicitly the
symmetry between time and frequency  than the classical conditions.
The result is first  stated and proved for $\Lz$, and then
generalized to coorbit spaces. As special cases, we obtain new
characterizations of compactness in  Besov-Triebel-Lizorkin spaces,
modulation spaces and  Bargmann-Fock spaces.
\end{abstract}

\section{Introduction}\label{int}

Compactness in function spaces is usually characterized by
conditions of the Arzela-Ascoli type. Typically, what is necessary
is an equicontinuity condition with respect to the norm of the
space under consideration. If the underlying topological space is
not compact, then in addition  a tightness condition is required,
i.e., all functions have the same ``essential''  support. The
prototype of such a  result  is the characterization of
compactness in $L^p$-spaces, which in its general form on locally
compact abelian groups is due to A. Weil. In the sequel $\chi _U$
will denote the indicator function of a compact set $U$.

\begin{Th}[\cite{Weil40}]\label{weil}
A closed and bounded subset $S$ of $\mathbf{L}^p (\mathbb{R}^d )$ for
$1\leq p  < \infty $
is compact if and only if the following conditions are satisfied:

(i) {\textsf  Equicontinuity:} for all $\epsilon > 0$ there
  exists $\delta > 0$ such that
  \begin{equation}
    \label{e1}
\sup_{f\in S} \sup _{|h| \leq \delta } \| f( . - h) - f \|_p
< \epsilon \, .
  \end{equation}

(ii)  {\textsf  Tightness:} for all $\epsilon > 0$ there exists a
compact set $U$ in $ \mathbb{R}^d$, such that
\begin{equation}
\label{t1}
\sup_{f\in S} \|f\,
\chi _U  -f\|_p < \epsilon \, .
\end{equation}
\end{Th}

  Far-reaching generalizations of Theorem~\ref{weil}  for general translation
invariant Banach spaces of distributions with a  so-called
double module structure    were proved in \cite{Feicompact84}.

Specializing to $\Lz$, it is well-known and not difficult to see
that the equicontinuity condition~\eqref{e1} is equivalent to the
tightness of the Fourier transforms $\widehat{S} = \{ \hat{f}: f
\in S \} $ in $\Lz$. In particular, \begin{it}
 a closed and bounded set  $S \subseteq \Lz$ is compact
\fif\ $S$ and $ \widehat{S}$ are both tight in $\Lz$, see
\cite{Feicompact84} and \cite{Pegocompact85}.\end{it}

The symmetry of this characterization under the \ft\ motivated us to
look at analytic tools which are designed expressedly to deal with
situations that treat a function and its \ft\ \emph{simultaneously} and
to search for a characterization of compactness by means of these tools.
In this regard the \stft\ is the tool that is used most frequently to
describe both time and frequency simultaneously, i.e., a function and
its \ft .

\begin{Def}[Short-time Fourier transform]
Let $M_{\omega }$ and $T_x$ denote frequency-shift by $\omega $
and time-shift by $x$, respectively, i.e., $ M_{\omega } T_x g(t) =
e^{2\pi i t\omega }g(t-x)$ for  $(x,\omega )\in \mathbb{R}^{2d}$.
The short-time Fourier transform (STFT)  of a function $f\in\mathbf{L}^2
(\mathbb{R}^d) $ with respect to a window function
$g\in\mathbf{L}^2 (\mathbb{R}^d) $ is defined as
\begin{equation}
  \label{ee2}
  \mathcal{S}_g f(x,\omega) = \int _{ \rd } f(t) \bar{g} (t-x) e^{-2\pi
  i \omega \cdot t} \, dt = \langle f,   M_{\omega } T_x g\rangle \,
  .
\end{equation}
\end{Def}
With slightly different normalization, the \stft\ also occurs
under the names ``(radar) ambiguity function'' or ``(cross-)Wigner
distribution'', see~\cite{book}.
 For suitable windows $g$, e.g. $g$ in the Schwartz class $
 \mathcal{S}(\rd )$, the value $   \mathcal{S}_g f(x,\omega) $
can be interpreted as a measure for the  energy of  $f$ at  $z =
\phas \in \mathbb{R}^{2d}$. An important property in the study of
compactness is the isometry property of the \stft\ which states
that for any $f,g\in\Lz$:
\begin{equation}
  \label{e3}
\|\mathcal{S}_g f\|_2 = \|g\|_2\|f\|_2 \, .
\end{equation}

It
becomes intuitively obvious that a condition comprising the
support conditions given in Theorem~\ref{weil} for time and
frequency separately can be formulated as a simultaneous tightness
condition in time and frequency via the short-time Fourier
transform.

\begin{Th}[Compactness in $\mathbf{L}^2 (\mathbb{R}^d
)$]\label{mainth} For a closed and bounded set $S\subseteq
\mathbf{L}^2 (\mathbb{R}^d )$  the following statements are
equivalent: \\
(i) S is compact in $\mathbf{L}^2 (\mathbb{R}^d )$. \\
(ii)  The set $\{ \mathcal{S} _g f : f \in S\} $ is tight in $\mathbf{L}^2
(\rdd ) $, this     means that for all $\epsilon > 0$ exists a compact
set $U\subseteq  \mathbb{R}^{2d}$, such that
\begin{equation}
\label{L2ST-tight}
\sup_{f\in S} \Big(\int_{U^c} |\mathcal{S}_g f(x,\omega) |^2 \, dx d\omega
\Big)^{\frac{1}{2}} < \epsilon .
\end{equation}
\end{Th}

\begin{proof}
To get an idea about possible generalizations we give the pretty proof
of this theorem right here. Without  loss of generality
we assume  that $ \|g\|_2 =1$, so that $\mathcal{S}_g$ is an isometry
on $\mathbf{L}^2 (\rd )$.

(i)$\Rightarrow$(ii): By compactness of $S$ we can find $f_1,
\ldots ,f_n$ such that \[\min_{j= 1,\ldots ,n} \| f-f_j \|_2 <
\frac{\epsilon}{2} \quad \quad \mbox{ for all } f\in S \, .\]
Since $\mathcal{S}_g f_j \in \mathbf{L}^2(\rdd )$ by \eqref{e3}, we
may choose a
compact set $U\subseteq \mathbb{R}^{2d}$ such that $ \int_{U^c}
|\mathcal{S}_g f_j (x,\omega) |^2 dx d\omega < \epsilon ^2  / 4 $
for $j = 1,\ldots ,n$. 
By \eqref{e3}  we obtain
 for arbitrary $f\in S$ that
\begin{eqnarray*}
 \lefteqn{\Big(\int_{U^c} |\mathcal{S}_g f (x,\omega) |^2 dx d\omega
   \Big)^{\frac{1}{2}} }\\
& \leq & \min_{j=1,\ldots , n}
   \Big\{ \Big( \int_{U^c} |\mathcal{S}_g (f-f_j ) (x,\omega) |^2 dx
   d\omega \Big)^{\frac{1}{2}} +   \Big(\int_{U^c} |\mathcal{S}_g f_j
   (x,\omega) |^2 dx d\omega \Big)^{\frac{1}{2}}\Big\}\\
   &\leq & \min_j\| f-f_j\|_2 +\frac{\epsilon}{2}<\epsilon
\end{eqnarray*}

(ii)$\Rightarrow$(i): It suffices to show  that every sequence $(f_n
)$ in $S$ contains a convergent subsequence.
By (ii) we can choose  a compact set $U \subseteq \mathbb{R} ^{2d}$  such
that
\begin{equation}\label{eq2}
\int_{U^c} |\mathcal{S}_g f(x,\omega) |^2 dx d\omega<\epsilon ^2
\end{equation}
holds for all $f\in S$, in particular for the sequence $(f_n)$.
Since by assumption $S$ is bounded, it is weakly compact in
$\mathbf{L}^2 (\mathbb{R}^d )$ and thus $(f_n )$ possesses a
weakly convergent subsequence $f_j = f_{n_j}$ with limit $f$,
i.e., $\langle f_j , h\rangle\rightarrow\langle f , h\rangle $ for
all $ h\in \mathbf{L}^2 (\mathbb{R}^d )$. Choosing $h=  M_{\omega
} T_x g$ for  $(x,\omega )\in \mathbb{R}^{2d}$, this implies the
pointwise convergence of the short-time Fourier transforms
\begin{equation}\label{conv}
\mathcal{S}_g f_j (x,\omega)\rightarrow\mathcal{S}_g f (x,\omega)
\quad \text{ for } x,\omega \in \rd \, .
\end{equation}
Since by \eqref{ee2} and the Cauchy-Schwarz inequality we have for
all $(x,\omega)$
\[|\mathcal{S}_g (f-f_j)(x,\omega)|\leq\|f-f_j\|_2\leq\sup_j\| f_j\|_2+\| f\|_2 < C ,\]
the restriction of  $|\mathcal{S}_g(f-f_j)|^2$ to $U$ is dominated
by the constant function $C^2 \chi _U\in \mathbf{L}^1
(\mathbb{R}^d)$. In view of (\ref{conv}) we may now  apply the
dominated convergence theorem and  obtain:
\begin{equation}\label{eq1}
\int_U |\mathcal{S}_g (f-f_j)(x,\omega)|^2 \,  dx
d\omega\rightarrow 0\ .
\end{equation}

The combination of  (\ref{eq1}) and   (\ref{eq2}) now  yields
\begin{eqnarray*}
\lefteqn{\overline{\lim} _{j\to \infty } \, \| f-f_j\|_2 =
  \overline{\lim} _{j\to \infty } \, \|\mathcal{S}_g
  (f-f_j) \|_2 }\\
  &\leq & \overline{\lim_{j\to \infty }}\Big( \int_U |\mathcal{S}_g (f
  -f_j)(x,\omega )|^2 dx d\omega \Big)^{\frac{1}{2}}+
  \overline{\lim_{j\to \infty }}\Big(\int_{U^c} |\mathcal{S}_g (f
  -f_j)(x,\omega )|^2 dx d\omega \Big)^{\frac{1}{2}}\\
  &\leq &  0+2\epsilon
\end{eqnarray*}
Therefore $\lim_{j\rightarrow\infty}\|f-f_j\|_2 = 0$
and thus $S$ is compact.
\end{proof}

Theorem~\ref{mainth} and its proof  suggest  several extensions. On
the one hand,  we may  replace  the
$\mathbf{L}^2$-norm of the short-time Fourier transform by other
norms and ask  for  which function spaces  we can still  characterize
compactness as in Theorem~\ref{mainth}. Pursuing this idea leads to
the characterization of compactness in the so-called modulation spaces
(Section~3.1).

On the other hand, if we are willing  to give up the
time-frequency interpretation of Theorems~\ref{weil} and~\ref{mainth},
we may  replace the \stft\  by other transforms.
As a  further  example occurring in modern analysis  we  consider
the wavelet transform, which shares the important  isometry property
with the STFT~\cite{Dau90}.
\begin{Def}[Continuous wavelet  transform]\label{CWT}
Let $T_x D_s  g(t) = s^{-\frac{d}{2}}g(s^{-1}(t-x))$ for
$(x,s)\in \mathbb{R}^{d}\times\mathbb{R}^+$. The continuous
wavelet transform of a function $f\in\mathbf{L}^2 (\mathbb{R}^d) $
with respect to a wavelet $g\in\mathbf{L}^2 (\mathbb{R}^d) $ is
defined to be
\begin{equation}
  \label{e24}
  \mathcal{W}_g f(x,s) = s^{-\frac{d}{2}} \int _{\rd } f(t)
  g(\frac{t-x}{s}) \, dt =    \langle f, T_x D_s g\rangle \, .
\end{equation}
\end{Def}
If $g$ is radial and satisfies the admissibility condition
\begin{equation}\label{wavcond}
\int _{\bR ^+ } |\hat{g}( t \omega )|^2\frac{dt }{t}  = 1  \quad \text{
  for all } \omega \in \rd \setminus \{0\} \, ,
\end{equation}
then $\int |\mathcal{W}_g f(x,s)|^2dx\frac{ds}{|s|^2}
  =\|f\|^2$ and thus 
  $\mathcal{W}_g$ is an isometry for $\Lz$~\cite{Dau90}.

The same proof as for Theorem~\ref{mainth} with the wavelet transform in place
of the STFT now yields the following criterion for compactness in $\lrd
$.

\begin{Th}[Compactness in $\mathbf{L}^2 (\mathbb{R}^d )$ via wavelet
  transform]
\label{waveth}
Let $g\in\Lz$ satisfy condition (\ref{wavcond}). A closed and
bounded set $S\subseteq \mathbf{L}^2 (\mathbb{R}^d )$ is compact
 in $\mathbf{L}^2 (\mathbb{R}^d )$ if and only if
for all $\epsilon > 0$ there exists a compact set $U\subseteq
\mathbb{R}^{d}\times\mathbb{R}^+$, such that
\begin{equation}
  \label{et6}
\sup_{f\in S} \Big(\int_{U^c} |\mathcal{W}_g f(x,s) |^2 \,
\frac{dxds}{s^{d+1}} \Big)^{\frac{1}{2}} < \epsilon \, .
\end{equation}
\end{Th}

The STFT and the wavelet transform do have other properties in common.
Both are defined as the inner product of $f$ with the action of a group
of unitary operators on a fixed function $g$. More precisely, both the
STFT and the wavelet transform are  \emph{ representation
  coefficients} of a certain  unitary continuous
 representations $\pi$ of a group $\mathcal{G}$ on a Hilbert space
 $\mathcal{H}$.
This observation has been very fruitful for the evolution of a general
wavelet theory~\cite{fg89jfa,fg89mh,GrochDF91}. In our context  we shall take
the proof of Theorem~\ref{mainth} as
an outline to obtain compactness criteria for a general class of
function spaces defined by means of other group representations.

To each irreducible unitary continuous representation $\pi $ of a locally
compact group on a Hilbert space $\cH $ satisfying some additional
integrability condition,  we associate a family of abstract function
spaces, the so-called \emph{coorbit spaces}.  For these spaces we will
prove compactness criteria analogous to those of Theorems~\ref{mainth} and
\ref{waveth}. Upon choosing a particular group and a natural representation,
we will recover the above  statements. In  addition we will
 obtain compactness criteria for modulation
spaces  by means of the short-time Fourier transform,
similar to Theorem~\ref{mainth},  and using the wavelet transform
we will  characterize compactness in
Besov-Triebel-Lizorkin spaces. Another modification leads to a new
compactness criterion for Bargmann-Fock spaces.

The paper is organized as follows.  Section~2  introduces the
concept of coorbit spaces and deals  with technical difficulties
arising in the generalization of Theorem~\ref{mainth} to coorbit
spaces. These concern the validity of dominated convergence and
the theorem of Alaoglu-Bourbaki. In Section~\ref{mainsec} we will
state the main theorem for general coorbit spaces. In
Section~\ref{Ex} we  treat the application of this theorem to
several  classes of well-known function spaces.

In a subsequent project we will apply the new compactness criteria to
study operators on coorbit spaces.

\section{Coorbit spaces} \label{coob}

 \subsection{Preliminaries and definition}\label{prel}

 We first recall the theory of coorbit spaces. For simplicity  we
 omit some technical details and  refer the reader
 to~\cite{fg89jfa,fg89mh,GrochDF91}  where the theory has been thoroughly
 investigated.

 The theory of coorbit spaces requires the  following basic structures:
 \begin{itemize}
 \item a locally compact  group $\mathcal{G}$   with Haar measure $dz$,
 \item   an irreducible  unitary   representation $\pi$ of $\mathcal{G}$
 on a Hilbert space $\mathcal{H}$,
 \item and a continuous  submultiplicative weight function
   $\nu$ on $\cG $, i.e., $\nu $  satisfies   $\nu (z_1 +z_2 )~\leq
   \nu (z_1 )\nu (z_2 )$ and $\nu(z_1) \geq 1$   for all $ z_1,  z_2\in\mathcal{G}$.
 \item A Banach space $(Y, \|. \|_Y)$ of functions on
   $\mathcal{G}$.
\end{itemize}

Functions of the form $z \in \cG \mapsto \langle f, \pi (z)
g\rangle$ are
 called  \emph{representation coefficients}  of $\pi$.
Upon inspection we see that the short-time Fourier transform
defined in~\eqref{ee2}  is (up to a trivial factor) a
representation
 coefficient of the Schr\"odinger representation of the Heisenberg
 group  $\mathcal{G}=\mathbb{R}^d\times\mathbb{R}^d\times\mathbb{T}$
 on $\Lz$,
 given by $\pi  (x,y,\tau ) f(t) = \tau e^{2\pi iy (t -x)}f(t -x)$.
 Likewise  the wavelet transform  is a
representation coefficient of the $ax+b$-group $\mathcal{G} =
\mathbb{R}^d\times\mathbb{R}^+ $  of the representation  $\rho (x,s)
f(t) = s^{-\frac{d}{2}}f(s^{-1}(t-x))$.

For  reasons of compatibility and well-definedness we impose the
following conditions on $\cG , (\pi  , \cH ), \nu $. We refer to
\cite{fg89jfa} for a detailed justification of the assumptions
stated above.

(A) $\pi $ is irreducible, unitary,  continuous  and
$\nu$-integrable, i.e., there exist   $g\in\mathcal{H}$, $g\neq
0$, such that
\begin{align}\label{propA}
   \int_{\mathcal{G}}|\langle\pi (z) g, g\rangle  |~\nu (z) \,  dz&
   <\infty \, .
\end{align}

(B) $Y$ is a \emph{solid} Banach function space  on $\cG $, i.e.,
if  $F\in Y$ and $G$ is measurable, satisfying
 $|G(z)| \leq |F(z)|$ for almost all $z\in
\mathcal{G}$, then  $G\in Y$ and $\|G\|_Y \leq \|F\|_Y$.

(C) $Y$ is invariant under right and left translations and
satisfies the convolution relation
$Y\ast\mathbf{L}^1_{\nu}(\mathcal{G})\subseteq Y$, with $\|F\ast
G\|_Y\leq\|F\|_Y\|G\|_{\mathbf{L}^1_{\nu}(\mathcal{G})}$ for $F\in
Y$, $G\in\mathbf{L}^1_{\nu}(\mathcal{G})$.

It follows  that $\mathbf{L}^{\infty}_0(\mathcal{G})$,  the space
of bounded function with compact support on $\mathcal{G}$ is
contained in $Y$. This property will be crucial in the proof of
our main statement.

We introduce the following notation for the representation
coefficient of $\pi$:
\[\mathcal{V}_g f(z)= \langle f, \pi (z) g\rangle  \mbox{ for }
z\in\mathcal{G}.\]

\begin{Def}[Abstract test functions and distributions] Fix
  $g_0\in \cH \setminus \{0\} $
  satisfying~\eqref{propA}.  Then the  space of test functions
$\mathcal{A}_{\nu}$ is defined as
\[\mathcal{A}_{\nu}
=\{g\in\mathcal{H}:\|g\|_{\mathcal{A}_{\nu}}=\|\mathcal{V}_{g_0}g
\|_{\mathbf{L}^1_{\nu}(\mathcal{G})}<\infty\}\]
\end{Def}
Then $\mathcal{A}_{\nu}$ is  dense in $\cH$, its
dual space  $\mathcal{A}_{\nu}'$, the space of all (conjugate-)
linear, continuous functionals on $\mathcal{A}_{\nu}$, contains $\cH$
and  plays the role of a space of distributions. It  will serve us as
a  reservoir of selection.
\begin{Def}[Coorbit spaces]\label{coodef}
Under the hypothesis (A), (B), (C) imposed  on $\cG , \pi, \cH , \nu
$, fix any $g
\in \cA _\nu \setminus \{0\}$. Then  the coorbit space of $Y$
under the representation $\pi$ is defined as
\[\mathcal{C}o_{\pi}Y\ = \{ f\in\mathcal{A}_{\nu}':~\RC  \in Y\}\]
with  norm 
$\|f\|_{\mathcal{C}o_{\pi}Y} = \|\RC \|_Y$.
\end{Def}

Then $\mathcal{C}o _\pi Y $ possesses the following properties,
see~\cite{fg89jfa}  for details.

(i) $\mathcal{C}o_{\pi}Y$ is a Banach space and invariant under
the action of $\pi $. Specifically,
\begin{equation}
  \label{ee?}
  \| \pi (z) f\|_{\mathcal{C}o _\pi Y} \leq C \nu (z)
  \|f\|_{\mathcal{C}o _\pi Y}    \quad \quad \text{ for } f \in
  \mathcal{C}o _\pi Y   \, .
\end{equation}

(ii) The definition of  $\mathcal{C}o_{\pi}Y $ is independent of the
 choice of $g\in \mathcal{A}_{\nu}$.

(iii)  Different functions  $g \in \cA _\nu \setminus \{0\} $ define
equivalent norms on $\mathcal{C}o_\pi Y$.

(iv) By definition $\mathcal{C}o_{\pi}Y$ is a subspace of $\cA
_\nu '$ and we also have
\begin{equation}
  \label{e6}
\|f \|_{\cA _\nu '} \leq C \| f \| _{\mathcal{C}o_{\pi}Y}
\end{equation}
(v)  Special cases:
$\mathcal{C}o_{\pi}\mathbf{L}^2(\mathcal{G})=\mathcal{H}$,
$\mathcal{C}o_{\pi}\mathbf{L}^1_{\nu}(\mathcal{G})=\mathcal{A}_{\nu}$,
and
$\mathcal{C}o_{\pi}\mathbf{L}^{\infty}_{\frac{1}{\nu}}(\mathcal{G})=\mathcal{A}_{\nu}'$.\\[.3cm]

In order to obtain compactness criteria analogous to
Theorem~\ref{mainth} for general coorbit spaces, we have to impose
further assumptions on $Y$. As pointed out at the end of
Section~\ref{int}, we need to use a norm $\| \cdot \|_Y$ for which
dominated convergence holds. In our treatment of dominated
convergence we follow~\cite[Ch.~1.3]{bennett}.

\begin{Def}
A Banach  function space   $Y$ on $\cG $ is said to have absolutely
continuous
norm, if $\| f\chi _{E_n}\|_Y\rightarrow 0$ for all $f$ and for
every sequence $\{ E_n\}_{n=1}^{\infty}$ of measurable subsets of
$\cG $ satisfying $E_n\rightarrow \emptyset$ almost
everywhere with respect to Haar measure.
\end{Def}
\begin{Pro}[\cite{bennett}] \label{abscont}
For a Banach  function space   $Y$  the following are equivalent.

(i) $Y$ has absolutely continuous norm.

(ii)  Dominated convergence holds for all $f~\in Y$: If  $f_n\in Y$,
$n=1,2,\ldots $,  and $g \in Y$ satisfy
$|f_n|\leq |g|$ for all $n$ and $f_n (z) \rightarrow f (z) $ a.~e.,
then $\| f_n -f \|_Y\rightarrow 0.$

(iii) The dual space $Y'$ of $Y$ coincides with its
associate space $Y^{\ast}$ defined as
\begin{equation}
  \label{e7}
  Y^{\ast}= \{g \mbox{ measurable, }
\sup_{f\in Y, \| f\|_Y \leq1}\int_{\cG}|f(z) g(z) |\, dz  <\infty\}
\end{equation}
\end{Pro}

\textbf{Example (Mixed-norm spaces):}
Let $m $ be a weight function on $\mathbb{R}^{2d}$ and let $1\leq
p,q \leq \infty$. Then the weighted mixed norm space $\Lmp$ consist
of all measurable functions on $\mathbb{R}^{2d}$ such that  the
norm
\[\|F\|_{\mathbf{L}^{p,q}_{m}}= \left (\int_{\mathbb{R}^d}\left (\int_{\mathbb{R}^d}
|F(x,\omega)|^p m (x, \omega )^p \, dx \right )^{\frac{q}{p}}\, d\omega \right )^
 {\frac{1}{q}}\] is finite, with the usual modifications when
$p=\infty $ or $q=\infty $.

If $m$ is a ``moderate'' weight with respect to the submultiplicative
weight $\nu $, i.e., $m(z_1+z_2) \leq C \nu (z_1 ) m(z_2)$, then
hypotheses (B) and (C) are always satisfied
(see~\cite[Prop. 11.1.3]{book}).

If $p,q < \infty$, then $L^{p,q}_{m} $ has an absolutely
continuous norm. For   $p=q=1$  this is just  Lebesgue's theorem
on dominated convergence. If $p,q < \infty$, then the dual space
is $L^{p',q'}_{1/m }$, where $1/p+1/p'=1$ \cite{BP}. As a
consequence of H\"older's inequality  the dual space coincides
with the associate space defined in  \eqref{e7}. By
Proposition~\ref{abscont} $L^{p,q}_{m}$ possesses an absolutely
continuous norm.

\subsection{Compactness in  coorbit spaces}\label{mainsec}

 We are now ready to state and prove our
main  theorem, a criterion for compactness in coorbit spaces.

\begin{Th}[Compactness in $\mathcal{C}o_{\pi}Y$]\label{coobth}
In addition to the general assumptions (A), (B), and (C),  assume that
$Y$ has an
absolutely continuous norm.  For a closed and bounded set
$S\subseteq \mathcal{C}o_{\pi}Y$  the following statements are
equivalent.

(i)  S is compact in $\mathcal{C}o_{\pi}Y$.

(ii) For all $~\epsilon> 0$ exists a compact set $U\subseteq
\mathcal{G}$, such that
\begin{equation}
  \label{e8}
  \sup_{f\in S}\|~\chi_{U^c}\, \RC \|_Y < \epsilon \, .
\end{equation}
\end{Th}

\begin{proof}
   The argument follows the simpler proof of Theorem~\ref{mainth}.

$(i)\Rightarrow(ii)$ Assume that $S$ is compact in
$\mathcal{C}o_\pi Y$ and let $\epsilon > 0$. Then there  exist
$f_1,\ldots, f_n \in S$ such that
$$
\min_{j=1, \dots , n}\|f-f_j\|_{\mathcal{C}o_{\pi}Y}
<\frac{\epsilon}{2} \quad \quad \text{ for all } f\in S \, .
$$
 Since   $\mathbf{L}^{\infty}_0 (\mathcal{G}) $ is
contained in $Y$ as a consequence of (B) and (C)  and since  $Y$ has
absolutely continuous norm,
$\mathbf{L}^{\infty}_0 (\mathcal{G}) $ is even dense in $Y$,
see~\cite[Prop.1.4]{Feicompact84}. Hence, there exist $H_j \in
\mathbf{L}^{\infty}_0\subseteq Y$ with $\|H_j - \mathcal{V}_g f_j
\|_Y <\frac{\epsilon}{2}$. By solidity of $Y$, $H_j $ can be
chosen as  the restriction $\chi_{U}\cdot \mathcal{V}_g f_j  $ for some
compact set $U\subseteq\mathcal{G}$,  and thus we obtain that
$$
\|\chi_{U^c}\mathcal{V}_g f_j\|_Y<\frac{\epsilon}{2}\quad \quad \text{
  for } j=1, \dots , n \, .
$$
 Then  for general  $f\in S$ we find that
\begin{eqnarray*}
\|\chi_{U^c}\mathcal{V}_g f  \|_Y &\leq &  \min_{j=1,\ldots ,n}
\Big( \|\chi_{U^c}\mathcal{V}_g (f-f_j) \|_Y + \|\chi_{U^c}\mathcal{V}_g
f_j\|_Y \Big) \\
& \leq &\min_{j=1,\ldots ,n}
\| \mathcal{V}_g (f-f_j) \|_Y + \frac{\epsilon}{2} \\
&=&  \min_{j=1,
  \dots , n} \|f- f_j\|_{\mathcal{C}o_{\pi}Y} + \frac{\epsilon}{2}  < \epsilon
\end{eqnarray*}
In  the second inequality  we have applied  condition (B) to the pointwise estimate
\[|\chi_{U^c}\mathcal{V}_g (f-f_j)  (z)  |\leq |\mathcal{V}_g (f-f_j)(z)
|\, .\]

$(ii)\Rightarrow(i)$
Assume that \eqref{e8} holds. Again it suffices to show that every
sequence $(f_n ) \subseteq S$ contains a convergent subsequence.

To extract a weak-star convergent subsequence of $(f_n) \subseteq S$, we
 modify the argument for $L^2$ follows:  Since $S$ is bounded and
closed in $\mathcal{C}o _{\pi }Y$, it is also bounded and closed in $\cA
_\nu ' $ by \eqref{e6}, and therefore  $S $ is weak-star compact in
$\mathcal{A}_{\nu}'$  by Alaoglu's theorem.
Consequently  we can find a weak-star convergent subsequence
$f_{n_j}$ of $(f_n)$, which we again denote by $f_j$,
with limit $f_\infty$ in $ S$, i.e. $\langle f_j , h \rangle
\longrightarrow \langle f_\infty , h \rangle  $  for all
$h~\in\mathcal{A}_{\nu}$.
In particular, for $h = \pi (z) g$, we obtain  \emph{pointwise}
convergence of  the representation coefficients on $\mathcal{G}$
$$
\mathcal{V}_g f_j (z) \longrightarrow \mathcal{V}_g f_{\infty} (z)
\quad \quad \text{for all } z \in \cG \, .
$$

Next we show that the sequence $\{ \mathcal{V}_g(f_{\infty} -f_j),
j \in \mathbb{N} \}$ is uniformly bounded on any compact set
$U\subseteq \cG $. We have
\[ |\langle  f_{\infty} -f_j , \pi (z) g\rangle  |\leq \|
f_{\infty} -f_j\|_{\mathcal{A}_{\nu}'}~\|\pi (z)
g\|_{\mathcal{A}_{\nu}},
\]
by duality, and  $\|\pi (z) g\|_{\mathcal{A}_{\nu}}\leq \nu
(z)\|g\|_{\mathcal{A}_{\nu}}$ by \eqref{ee?}. Therefore
$$
\sup _{z\in U} |\mathcal{V}_g(f_{\infty} -f_j)(z)| \leq \|g\|_{\cA
_\nu } \sup _{z\in U} \nu (z) \,  \sup _{j\in \bN } \|f_{\infty}
- f_j\| _{\cA _\nu '} \le C \chi _U(z)
$$
Since $\chi _U \in Y$ and $Y$ has an absolutely continuous
norm, we can apply dominated convergence (Proposition~\ref{abscont} ii)) to
obtain
\begin{equation}
  \label{e23}
\lim _{j\to \infty  } \|\chi _U \cdot  \mathcal{V}_g (f_{\infty}
-f_j)\|_Y = 0 \, .
\end{equation}

To deal with the behavior of $\mathcal{V}_g(f_{\infty} -f_j)$ on
the complement $U^c$, we use the assumption \eqref{e8}. Given
$\epsilon >0$, we choose $U \subseteq \cG $ so that $ \|\chi
_{U^c} \, \mathcal{V} _g f \|_Y < \epsilon /2$ for all $f \in
S\cup  \{ f_{\infty} \}$.
The combination of these steps now yields
\begin{align*}
  \lim_{j \to \infty } \| f_{\infty} -f_j\|_{\mathcal{C}oY}&=\lim_{j \to \infty
    }\| \cV _g( f_{\infty} -f_j)  \|_Y\\
  &\leq \lim_{j \to \infty } \| \chi _U\cdot  \cV _g( f_{\infty} -f_j)  \|_Y+
  \overline{\lim } _{j     \to \infty }\| \chi _{U^c} \cdot \cV _g( f_{\infty} -f_j)
  \|_Y  \\
  &\leq 0+2 \sup _{f \in S\cup  \{f_{\infty}\}} \| \chi _{U^c} \cdot \cV _g f  \|_Y
<  2\epsilon \, .
\end{align*}
Therefore any sequence in $S$ has a  subsequence that converges in
$\mathcal{C}o_{\pi}Y$  and so $S$ is compact.
\end{proof}

\noindent \textbf{Remarks:} 1. Loosely speaking, Theorem~\ref{coobth}
states that  a set in
$\mathcal{C}o _\pi Y$ is compact, if and only if the set of
representation coefficients is tight in $Y$.

2. Note that in the first part of the proof   we have  only  used
 the fact
that  $\mathbf{L}^{\infty}_0(\mathcal{G})$ is dense in $ Y$.  On
the
 other hand, the absolutely continuous
norm  of $Y$ is only needed for the proof of sufficiency of
{\it(ii)} for $S$ to be compact.
If $\mathbf{L}^{\infty}_0(\mathcal{G})$ is not dense in $Y$, then
condition {\it(ii)} characterizes the compactness in the closed
subspace $\mathcal{C}o_{\pi} Y_0$  of $\mathcal{C}o_{\pi} Y$,
where $ Y_0$
 is the closure of
$\mathbf{L}^{\infty}_0(\mathcal{G})$ in $Y$.

3.  If  $Y$ does not possess an absolutely continuous norm, then
one may alternatively apply
the compactness criterion  to the coorbit corresponding to  the closed
subspace $Y_a\subseteq  Y$ of all ``functions of absolutely
continuous norm'' (see \cite[Ch.~1.3]{bennett}).

\section{Examples}\label{Ex}

Theorem~\ref{coobth} yields a handy compactness criterion for most
function spaces commonly used in analysis. We now give several
concrete manifestations of Theorem~\ref{coobth}.
In order to apply it
we have to verify  that all the conditions on $\cG ,
\pi , \nu $ and $Y$ are  fulfilled. Once the general setting is
described, this is  an easy task.

\subsection{Modulation spaces}
Modulation spaces are those function spaces which are associated to the
short-time Fourier transform~\eqref{ee2}.

Their standard definition is as follows. Fix a non-zero "window
function" $g \in \cS (\rd )$   and consider moderate  
functions $m $ satisfying $m (z_1+z_2) \le C  (1+|z_1|)^s m
(z_2)$, $z_1,z_2 \in \rdd $ for some constants $C,s \ge 0$, for
instance $m(z) = (1+|z|)^a $ for $a \in \mathbb{R} $ is moderate
with respect to $\nu (z) = (1+|z|)^{|a|}$.  Then  the modulation
space $\MS$ is defined as the space of all tempered distributions
$f\in\SWD$ with $\mathcal{V}_g f\in \Lmp$, with  norm
 \[\| f \|_{\MS} = \|\mathcal{S}_g f\|_{\Lmp}.\]

For the detailed theory of the \modsp s we refer to \cite[Ch.
11--13]{book} where they are treated for even more general classes
of weight functions. As particularly important modulation space we
mention $M^{1,1}_m $ with constant weight $m \equiv 1$. In the
abstract notation it is just $\cA _\nu $, it is a Segal algebra
and is denoted by  $S_0$ in harmonic analysis.

To interpret \modsp s as coorbit spaces, we extend the \tfs s
$\phas \to T_x M_\omega $ to a unitary \rep\ of the Heisenberg
group. Let  $\mathbb{H}=
\mathbb{R}^d\times\mathbb{R}^d\times\mathbb{T}$ be the $d$-dimensional
reduced Heisenberg
group with multiplication \[ (x_1,\omega _1, e^{2\pi i\tau
_1})(x_2,\omega _2, e^{2\pi i\tau _2}) = (x_1 +x_2, \omega _1
+\omega _2,e^{2\pi i(\tau _1 +\tau _2)}e^{\pi i x_2\cdot \omega _1
})\] and let
 $\pi $ be the Schr\"odinger \rep\ of $\bH $  acting
on $\lrd $ by time-frequency-shifts
\begin{equation}
  \label{e26}
  \pi (x,\omega,\tau) =  e^{2\pi i\tau}T_x
M_{\omega} \, .
\end{equation}
Then $\pi $ is an irreducible, unitary representation of $\mathbb{H}$
on $\Lz$.
The \rep\ coefficient for the Gaussian  $\phi (t)= e^{-\pi t \cdot t}$
is $\langle \phi ,\pi (x, \omega , \tau )\phi  \rangle = 2^{-d/2} \, \bar{\tau
} e^{\pi i x \cdot \om } \, e^{-\pi (x\cdot x + \om \cdot \om )}  $,
therefore $\pi $ is integrable with respect to any weight $\nu (x) =
\mathcal{O} (e^{\alpha |z|}), \alpha \geq 0$, see~\cite{book}.
Furthermore observe that
\[|\RC (x,\omega )| = |\langle f, \pi (x,\omega,\tau) g\rangle |.\]
Now consider  the auxiliary space $\widetilde{\mathbf{L}^{p,q}_{m}}$
consisting  of all measurable functions $f$  on $\mathbb{H}$
  such that
  \[\Big(\int_{\tau\in \mathbb{T}}|f(x,\omega ,\tau )|^2d\tau \Big)^{\frac{1}{2}}\in \mathbf{L}^{p,q}_{m} ,\]
then the modulation spaces can be interpreted as coorbit spaces by
\[\MS= \mathcal{C}o_{\pi}\widetilde{\mathbf{L}^{p,q}_{m}}  .\]

If $1\le p,q < \infty$, then  $\mathbf{L}^{\infty}_0 (\mathbb{H})$
is dense in $\widetilde{\mathbf{L}^{p,q}_{m}}$.  We have already
verified after Proposition~\ref{abscont} that the spaces
$\widetilde{\mathbf{L}^{p,q}_{m}}$ possess an absolutely continuous
norm. Therefore all conditions of Theorem~\ref{coobth} are satisfied, and we
obtain the following more explicit characterization of compactness in
\modsp s.

\newpage
\begin{Th}[Compactness in $\MS$]\label{modth}

Let $0 \neq g\in \mathbf{M}^1_{\nu } (\mathbb{R}^{d})$, $1\leq p,q <
\infty $ and $S$ be   a closed and bounded subset of $ \MS$.
Then $S$ is compact in $\MS$ if and only if for all $~\epsilon> 0$
exists a compact set $U\subseteq \rdd $, such that
\[\sup_{f\in S}\|~\chi_{U^c}\cdot \mathcal{V}_g f\|_{\mathbf{L}^{p,q}_{m}} < \epsilon .\]
\end{Th}

\textbf{Remark:} Clearly  such a characterization  cannot hold when
$p=\infty $ or
$q=\infty $. In this case $\MS $ is the dual of a non-reflexive Banach
space, and compactness in norm takes on a different shape.

\subsection{Besov-Triebel-Lizorkin spaces}
The Besov-Triebel-Lizorkin spaces
 are those function spaces that can be  associated to the wavelet
transform~\eqref{e24}. Let $g \in \cS (\rd )$ be a fixed nonzero
radial function with all moments vanishing. Then the homogeneous
Besov space $\BTL$ contains all tempered distributions (modulo
polynomials) such that
\begin{equation}
  \label{e25}
  \|f\|^q_{\btl} = \int_{\mathbb{R}}
\left (\int_{\mathbb{R}^d} |\langle f, \rho (x,s) g\rangle|^p
dx \right )^{\frac{q}{p}}
s^{-q(\alpha +\frac{d}{2}-\frac{d}{q})}\frac{ds}{s^{d+1}}<\infty
\end{equation}
This definition is equivalent to the standard definition given in
\cite{triebel83}, see \cite{triebel88}.

To interpret the Besov spaces as coorbit spaces, we consider the
$ax+b$--group $\cG = \mathbb{R}^d\times\mathbb{R}^+ $ with
multiplication $(b,a) \cdot (x,y) = (ax+b, ay), b,x \in \rd, a,y
\in \bR ^+,  $ and the  representation of $\cG $ on $\lrd $
 by  translations and dilations
\begin{align}\label{btlrep}
\rho  (x,s) f(t) = s^{-\frac{d}{2}}f(s^{-1}(t-x)).
\end{align}

Again it is easy to see that this representation is integrable
with respect to all weights of the form $\nu (x,s)  = \max
(1,s^\alpha ) $ for some
$\alpha \geq 0$ 
by choosing $g \in \cS (\rd )$ such that supp$\, \hat{g} \subseteq
\{t\in \rd: 0 < c \le |t| \le d  <  \infty \}       $. However, this \rep\
is reducible, and thus the general theory of Section~\ref{coob} is not
immediately applicable. To save the situation, we take the extended
group $\mathbb{R}^d\times (\mathbb{R}^+\times SO(d))$ with the \rep\
$\pi  (x,s,\cO )f(t)=
s^{-\frac{d}{2}}f(s^{-1}(\cO ^{-1}   (t-x))), \cO \in SO(d),$ acting on
$\Lz$. Then $\pi $ is again irreducible. Now
take a wavelet $g$ that is rotation invariant, then
\[\langle f,\pi  (x,s, \cO  )g\rangle=\langle f,\rho (x,s)g\rangle\]

Comparing with~\eqref{e25}, we see that
$$
\BTL = \cC o _\pi \mathbf{L}^{p,q}_{\alpha  +\frac{d}{2}-\frac{d}{q}}
$$
where the subscript refers to the weight $\nu (x,t, \cO ) =
t^{-(\alpha +\frac{d}{2}-\frac{d}{q})}$ on the extended $ax+b$--group.
As before all assumptions of Theorem~\ref{coobth} are satisfied, and
we obtain
the following new characterization of compactness in Besov spaces.

\begin{Th}[Compactness in $\BTL$]\label{btlth}
A closed and bounded set $S\subseteq \BTL$, $1\leq p,q <\infty$,
is compact in $\BTL$ if and only if for all $~\epsilon> 0$ there  exists
a compact set $U\subseteq \mathbb{R}^d\times\mathbb{R}^+ $, such
that
\[\sup_{f\in S}\|~\chi_{U^c}\cdot \mathcal{W}_g f\|_{\mathbf{L}^{p,q}_{\nu}} < \epsilon .\]
\end{Th}

Similarly, all Triebel-Lizorkin spaces $\TTL$, among them $L^p$
and the Hardy spaces,  can be defined as the coorbits of so-called
tent spaces $\mathbf{T}^{p,q}_{\nu}$ on $\mathcal{G}$,
cf.~\cite{CMS}. The  compactness in  $\TTL$ can be characterized
as in Theorem~\ref{btlth} with $\mathbf{L}^{p,q}_{\nu}$ replaced
by $\mathbf{T}^{p,q}_{\nu}$. Since the classical criterion of
Theorem~\ref{weil} is much simpler to use, we omit the explicit
formulation of Theorem~\ref{btlth} for the $\TTL $-spaces.


\subsection{Bargmann-Fock spaces}

Finally we study a class of function spaces occurring in complex
analysis, see \cite{JPR,folland89}.

\begin{Def}
The Bargmann-Fock spaces  $ \cF ^p = \BF$, $p<\infty$, are the Banach spaces
of entire functions $F$ on $\mathbb{C}^d$ for which the norm
\[\|F\|_{\BF}= \Big(\int_{\mathbb{C}^d}|F(z)|^p e^{-\frac{\pi}{2}
p|z|^2}dz\Big)^{\frac{1}{p}}, \] is finite.
\end{Def}
\textbf{Remark:} $\cF ^2 $ is a Hilbert space with inner product
\[\langle F,G\rangle _{\cF ^2 } = \int_{\mathbb{C}^d}F(z)
\overline{G(z)}e^{-\pi |z|^2}dz,\] which is isometrically
isomorphic to $\Lz$ via the Bargmann transform,
see~\cite[Ch.~3]{book}.

By identifying $\mathbb{H} = \mathbb{R}^{2d}\times\mathbb{T}$
with $\mathbb{C}^d\times\mathbb{T}$  and using the notation of
\cite[p.~183]{book}, the Heisenberg group  acts on $\cF ^2 $ via the
Bargmann-Fock representation $\beta $ as follows:
\[\beta (z, \tau ) F(w ) = e^{2\pi i\tau} e^{\pi z \cdot w}
F(w -\overline{z}) e^{-\frac{\pi |z|^2}{2}} \quad \quad z,w \in
\mathbb{C}^d, |\tau | = 1 \, .\]
Then $\beta $ is irreducible on $\cF ^2$ and is in fact equivalent to
the Schr\"odinger \rep\ $\pi $ of \eqref{e26}. Therefore $\beta $ enjoys
all properties  required to apply Theorem~\ref{coobth}. The following
compactness criterion for Bargmann-Fock spaces seems to be new.

\begin{Th}[Compactness in $\BF$]\label{bfth}
Let $1\leq p<\infty$.  A closed and bounded set $S\subseteq \BF$
 is compact in $\BF$
if and only if  for all $\epsilon> 0$  there exists a compact set
$U\subseteq \mathbf{C}^d$, such that
\[\sup_{F\in S}\Big(\int_{U^c}|F(z)|^p e^{-\frac{\pi}{2}
p|z|^2}dz\Big)^{\frac{1}{p}} < \epsilon .\]
\end{Th}

\begin{proof}
We show the identification   $\BF=\mathcal{C}o_{\beta}\Lp$ by using
the properties of $\cF ^2$ as a reproducing kernel Hilbert
space~\cite[Thm.~3.4.2]{book}.

 We know that $F(\xi ) =  \langle F,e^{\pi z \overline{\xi}}\rangle  $.
In particular, for the constant $\mathbf{1}$ we  obtain that
$$
|\langle \mathbf{1},\beta (z,\tau ) \mathbf{1}\rangle  |= |\langle
\mathbf{1},e^{\pi\xi \overline{z}}\rangle  e^{-\frac{\pi}{2}|z|^2}|=
e^{-\frac{\pi}{2}|z|^2}\, ,
$$
which implies that $\beta $ is integrable with respect to arbitrary
weights  $\nu (x) =
\mathcal{O} (e^{\alpha |z|})$ and that $\mathbf{1} \in \cA _\nu $.
Furthermore, the identity
\[|\langle F, \beta (z,\tau )\mathbf{1}\rangle  | = |\langle F,e^{\pi\xi
\overline{z}}\, e^{-\frac{\pi}{2}|z|^2}\rangle  | = |F(z
)|e^{-\frac{\pi}{2}|z|^2}\]
implies that
\[\langle F, \beta (z,\tau )\mathbf{1}\rangle  \in\Lp\Longleftrightarrow
F\in\BF .\]
This means that  $\BF = \mathcal{C}o_{\beta}\Lp$. Thus $\BF$  is a
coorbit space and  so  the statement follows from Theorem \ref{coobth}.
\end{proof}

\textbf{Remark:}
We leave it to the reader to   generalize the result  to \emph{weighted}
Bargmann-Fock spaces $\mathcal{F}^p_{\nu} (\mathbb{C}^{d} )$ or to
``mixed-norm'' Bargmann-Fock spaces.

\bibliographystyle{plain}
\bibliography{monibib,long_names,all_refs}
\end{document}

%% file: monimacro.tex
\newcommand{\Lz}{\mathbf{L}^2 (\mathbb{R}^d )}

\newcommand{\Lp}{\mathbf{L}^p (\mathbb{R}^d )}
\newcommand{\Lmp}{\mathbf{L}^{p,q}_{\nu} (\mathbb{R}^{2d} )}
\newcommand{\MS}{\mathbf{M}^{p,q}_m (\mathbb{R}^{d})} 
\newcommand{\BF}{\mathcal{F}^p (\mathbb{C}^d )}

\newcommand{\SWD}{\mathcal{S}'(\mathbb{R}^d )}
\newcommand{\BTL}{\dot{\mathbf{B}}_{p,q}^s (\mathbb{R}^{d})} 
\newcommand{\btl}{\mathbf{B}_{p,q}^s } 
\newcommand{\TTL}{\mathbf{F}_{p,q}^s (\mathbb{R}^{d})} 

\newcommand{\RC}{\mathcal{V}_g f}

\def\rd{\bR^d}
\def\phas{(x,\omega )}
\def\rdd{{\bR^{2d}}}
\def\lrd{\mathbf{L}^2 (\mathbb{R}^d )}

\newcommand{\field}[1]{\mathbb{#1}}
\newcommand{\bR}{\field{R}} 
\newcommand{\bN}{\field{N}} 
\newcommand{\bH}{\field{H}} %

\newcommand{\ft}{Fourier transform}
\newcommand{\stft}{short-time Fourier transform}

\newcommand{\fif}{if and only if}
\newcommand{\tfs}{time-frequency shift}

\newcommand{\modsp}{modulation space}

\newcommand{\rep}{representation}

\newcommand{\om}{\omega}

\def\cF{\mathcal{F}} 
\def\cS{\mathcal{S}}

\def\cH{\mathcal{H}}

\def\cG{\mathcal{G}}

\def\cA{\mathcal{A}}

\def\cC{\mathcal{C}}

\def\cO{\mathcal{O}}

\def\cV{\mathcal{V}}